\documentclass[12pt]{article}
\usepackage{mathptmx}
\usepackage{graphicx}
\usepackage{amsmath}
\usepackage{amsfonts}
\newtheorem{thm}{Theorem}[section]

\newtheorem{rmk}{Remark}[section]

\newtheorem{tb}{Table}
\usepackage{epstopdf}

\title{On Bivariate Kantorovich Exponential Sampling Series }
\author{ Prashant Kumar \thanks{Department of Mathematics, Visvesvaraya National Institute of Technology, Nagpur, Nagpur-440010, India. \newline E-mail: pranwd92@gmail.com}
 \and
 A. Sathish Kumar  \thanks{Department of Mathematics, Visvesvaraya National Institute of Technology, Nagpur, Nagpur-440010, India. \newline E-mail: mathsatish9@gmail.com}
  \and
  Shivam Bajpeyi \thanks{Department of Mathematics, Visvesvaraya National Institute of Technology, Nagpur, Nagpur-440010, India.
  \newline E-mail: shivambajpai1010@gmail.com}
  }
\date{}

\begin{document}
\maketitle
\bibliographystyle{plain}
\abstract{

We analyse the approximation properties of the bivariate generalization of the family of Kantorovich type exponential sampling series. We derive the point-wise and Voronovskaya type theorem for these sampling type series. Using the modulus of smoothness, we obtain the quantitative estimate of order of convergence of these series. Further, we establish the degree of approximation for these series associated with generalized Boolean sum (GBS) operators. Finally, we provide a few examples of kernels to which the theory can be applied along with the graphical representation and error estimates.

\endabstract



\section{Introduction}

The exponential sampling formula has been initiated by Bertero, Pike \cite{bertero} and Gori \cite{gori} and they establish a series representation for the class of Mellin band-limited function exploiting its exponentially spaced samples.
For $f:\mathbb{R}^+ \rightarrow \mathbb{C}$ and $c \in \mathbb{R},$ the exponential sampling formula is given by
$$ (E_{c,T}f)(x):= \sum_{k=-\infty}^{\infty} lin_{\frac{c}{T}}(e^{-k}x^{T}) f(e^{\frac{k}{T}})$$
where $lin_{c}(x)= \frac{x^{-c}}{2\pi i} \frac{x^{\pi i}-x^{-\pi i}}{\log c} = x^{-c} sinc(\log x)$ with continuous extension $lin_{c}(1)=1.$ The above sampling formula provides a very useful tool to approximate a signal by using its values at the node $ (e^{\frac{k}{w}})$ and this also called the Mellin-version of the \textit{Shannon sampling theorem} (see \cite{butzer2}).
Butzer and Jansche \cite{butzer5} were first initiated the mathematical study of exponential sampling formula exploiting the Mellin analysis. An independent study of Mellin theory was studied by Mamedov in \cite{mamedeo} and then developed by Butzer et.al. in \cite{butzer3,butzer4,butzer5,butzer7}. Subsequently, various authors have contributed in the advancement of the Mellin theory, see \cite{bardaro1,bardaro9,bardaro2,bardaro3}. The convergence of the generalized exponential sampling series and their linear combination was studied extensively in \cite{bardaro7,comboexp} respectively. Subsequently, the approximation properties of above series was studied in Mellin-Lebesgue spaces in \cite{bardaro11}. The approximation properties of the bivariate Kantorovich exponential sampling operators has been investigated in \cite{bvexp}. In order to reduce the time jitter-error, the Kantorovich version of the generalized sampling series was introduced and studied in \cite{own}. Let $ x \in \mathbb{R}^{+}$ and $ w >0.$ Then the Kantorovich type exponential sampling series is defined as
\begin{equation*}
(I_{w}^{\chi}f)(x):= \sum_{k= - \infty}^{+\infty} \chi(e^{-k} x^{w})\  w \int_{\frac{k}{w}}^{\frac{k+1}{w}} f(e^{u})\  du
\end{equation*}
where $ f: \mathbb{R}^{+} \rightarrow \mathbb{R}$ is locally integrable such that the above series is absolutely convergent for every $ x \in \mathbb{R}^{+}.$ This provides an useful tool to approximate approximate Lebesgue integrable functions by using its samples at the nodes $ (e^{\frac{k}{w}})_{w>0},\ k \in \mathbb{Z}.$ In the spirit of improving the order of approximation for the above family, the theory of linear combination of these operators was developed in \cite{lc} using the tools of Mellin analysis. The Kantorovich type modification of sampling series has been a topic of significant applications in approximation theory over the past decades. The Kantorovich modifications of several operators have been studied in various settings, see eg. \cite{bardaro10,bardaro6,costa2,costa3,vinti1,Dani,gupta2,deva}. The approximation behaviour of multivariate and non-linear Kantorovich type sampling operators was investigated in \cite{bv1,butz1,butz2,COS,Dan}. This paper deals with the bivariate extension of the Kantorovich exponential sampling series along with the associated GBS operators.\\
The paper is organized as follows. We introduce some basic definitions and preliminary results in Section 2. The study of convergence results, eg. point-wise convergence theorem and Voronovsakaya type theorem and its quantitative estimates in terms of modulus of continuity for the family (\ref{main}) has been given in Section 3. Subsequently, in Section 4, we define the GBS operators of bivariate Kantorovich exponential sampling series and establish some approximation results using the notion of Mellin B\"{o}gel continuous and Mellin B\"{o}gel differentiable functions. We also provide a few examples based on the presented theory along with graphical representation and error-estimates.

\section{Preliminaries}
Let $\mathbb{R}_{+}$ denotes the set of all positive real numbers and $\mathbb{R}^{2}_{+} := \mathbb{R}_{+} \times \mathbb{R}_{+}.$ Let $ C^{(r)}(\mathbb{R} ^{2}_{+})$ denotes the space of all continuous functions with upto $r^{th}$ order partial derivatives are continuous and bounded on $ \mathbb{R}^{2}_{+}$ equipped with the supremum norm $\|f\|_{\infty} := \sup_{(x,y) \in \mathbb{R} ^{2}_{+}} |f(x,y)|.$ For convenience, we set $C(\mathbb{R} ^{2}_{+}):= C^{(0)}(\mathbb{R} ^{2}_{+}).$ For any $\delta >0$ and $(x,y) \in \mathbb{R} ^{2}_{+},$ $B_{\delta}(x,y)$ denotes  the open ball of radius $\delta$ and centred at $(x,y),$ namely
$$B_{\delta}(x,y) := \{ (u,v) \in \mathbb{R} ^{2}_{+} : (x-u)^{2}+(y-v)^{2} < \delta^{2}\}.$$
A function $f: \mathbb{R}^{2}_{+} \rightarrow \mathbb{C}^{2}$ is called log-uniformly continuous on $\mathbb{R}^{2}_{+}$ if for any given
$\epsilon > 0$ there exists $\delta > 0$ such that $|f(x,y) -f(u,v)| < \epsilon $ whenever $ \left( (\log x - \log u)^{2}+(\log y - \log v)^{2} \right) < \delta^{2},$ for any $(x,y),\ (u,v) \in \mathbb{R}^{2}_{+}.$ We denote $\mathcal{C}(\mathbb{R}^{2}_{+})$ containing all log-uniformly continuous and bounded functions defined on  $\mathbb{R}^{2}_{+}.$ Moreover, $L^{\infty}(\mathbb{R}^{2}_{+}) $ denotes the space of all bounded functions on  $\mathbb{R}^{2}_{+}.$\\
The notion of Mellin partial derivatives of any function $f:\mathbb{R}^{2}_{+} \rightarrow \mathbb{C}$ with respect to variables $x,y \ $ is defined as
$$ \theta_{x}f := x \frac{\partial f}{\partial x} \ \ \ \ \mbox{and} \ \ \ \theta_{y}f := y \frac{\partial f}{\partial y}.$$
For $h=(h_{1},h_{2}) \in \mathbb{N}_{0}^{2},$ the Mellin partial derivatives of order $r,$ where $r=|h|=h_{1}+h_{2}$ is given by
$$ \theta^{r}_{x^{h_{1}}y^{h_{2}}}f := \theta_{x}^{h_{1}}(\theta_{y}^{h_{2}} f).$$ In particular, for $r=2$ we have
$\theta^{2}_{x}f:= \theta_{x}(\theta_{x}f)$ and $\theta^{2}_{y}f:= \theta_{y}(\theta_{y}f).$
Let $f:\mathbb{R}^{2}_{+} \rightarrow \mathbb{C}$ be such that $f \in C^{(r)}(\mathbb{R}^{2}_{+}),\ r \in \mathbb{N}.$ For $(x,y), (s,t) \in \mathbb{R}^{2}_{+},$ the bivariate Taylor formula in the Mellin setting (see \cite{bvexp}) is given by
$$ f(sx,ty)=f(x,y)+(\theta_{x} \log s+\theta_{y} \log t)f(x,y)+  \cdots + \frac{(\theta_{x} \log s+\theta_{y} \log t))^{r-1}}{(r-1)!}f(x,y)+ R_{r}(s,t),$$ where $R_{r}(s,t):= h(s,t) (\log^{2}s+ \log^{2}t)^{\frac{r}{2}}$ and $\displaystyle \lim_{(s,t)\rightarrow (1,1)} h(s,t) = 0.$

Assume that $\chi\in C( \mathbb{R}^{2}_{+})$ be fixed. Then for any $ \eta \in \mathbb{N}_{0}= \mathbb{N} \cup \{0 \},$ $ p= (p_{1},p_{2})\in \mathbb{N}_{0}^2$ with $ \vert p \vert = p_{1}+p_{2}= \eta,$ we define the algebraic moments of order $\eta$ as
\begin{eqnarray*}
m_{( p_{1},p_{2})} (\chi,u,v):=\sum_{k=-\infty}^{\infty}\sum_{j=-\infty}^{\infty}
\chi(e^{-k}u ,e^{-j}v)(k-\log u)^{p_{1}}(j-\log v)^{p_{2}}
\end{eqnarray*}
and the absolute moments by
\begin{eqnarray*}
M_{( p_{1},p_{2})} (\chi):=\sum_{k=-\infty}^{\infty}\sum_{j=-\infty}^{\infty}|\chi(e^{-k}u ,e^{-j}v)| |k-\log u|^{p_{1}}|j-\log v|^{p_{2}}
\end{eqnarray*}
Also,
$$ M_{\eta}(\chi):= \max_{\vert p \vert = \eta} \  M_{(p_{1},p_{2})} (\chi).$$

\begin{rmk}
We can easily see that for $ \xi,\eta \in \mathbb{N}_{0}$ with $\xi < \eta,$  $ M_{\eta}(\chi)< + \infty$ implies that $M_{\xi}(\chi)< + \infty.$ Indeed, for $p_{1}+p_{2}= \xi$

\noindent $\displaystyle \sum_{k=-\infty}^{\infty}\sum_{j=-\infty}^{\infty}|\chi(e^{-k}u ,e^{-j}v)|\, |k- \log u|^{p_{1}}\, |j-\log v|^{p_{2}}$
\begin{eqnarray*}
 &=& \sum \sum_{(k,j) \in B_{1}(\log u,\log v)} |\chi(e^{-k}u ,e^{-j}v)|\, |k-\log u|^{p_{1}}\, |j-\log v|^{p_{2}} \\&&
 +\sum \sum_{(k,j) \notin B_{1}(\log u,\log v)}|\chi(e^{-k}u ,e^{-j}v)|\, |k-\log u|^{p_{1}}\, |j-\log v|^{p_{2}} \\
&\leq & ( \|\chi \|_{\infty} + M_{\eta}(\chi)) \hspace{0.3cm} < + \infty.
\end{eqnarray*}
Also note that when $\chi$ is compactly supported then we have $ M_{\eta}(\chi)< \infty,$ for every $ \eta \in \mathbb{N}_{0}.$
\end{rmk}
Now we assume that the kernel satisfies the following conditions :
\begin{itemize}
\item[K1)] The series $ \displaystyle{\sum_{k=-\infty}^{\infty}\sum_{j=-\infty}^{\infty}} \chi(e^{-k} u, e^{-j}v) =1 ,$ for every $(u,v) \in \mathbb{R}^{2}_{+}.$

\item[K2)]$ M_{2}(\chi) < +\infty$ and
$$ \lim_{\gamma \rightarrow + \infty} \sum \sum_{(k,j) \notin  B_{\gamma}(\log u,\log v)}|\chi(e^{-k} u, e^{-j}v)| \ |k- \log u|^{p_{1}} \ |j- \log v|^{p_{2}}=0$$ uniformly for $ (u,v) \in \mathbb{R}^{2}_{+},$ where $p_{1}+p_{2}= 2.$
\end{itemize}

\begin{rmk}
The condition (K2) implies that
$$\lim_{\gamma \rightarrow + \infty} \sum \sum_{(k,j) \notin  B_{\gamma}(\log u,\log v)} |\chi(e^{-k} u, e^{-j}v)| \ |k- \log u|^{p_{1}} \ |j- \log v|^{p_{2}}=0 \ \ \mbox{for}\ \ \ p_{1}+p_{2}=h,\ h=0,1.$$
We easily see that

\noindent $\displaystyle \lim_{\gamma \rightarrow + \infty} \sum \sum_{(k,j) \notin  B_{\gamma}(\log u,\log v)}|\chi(e^{-k} u, e^{-j}v)| \ |k- \log u|^{p_{1}} \ |j- \log v|^{p_{2}}$
\begin{eqnarray*}
&\leq &  \lim_{\gamma \rightarrow + \infty} \frac{1}{\gamma^{2-h}} \sum \sum_{(k,j) \notin  B_{\gamma}(\log u,\log v)}|\chi(e^{-k} u, e^{-j}v)| \ |k- \log u|^{(p_{1}-h+2)} \ |j- \log v|^{p_{2}} \\
&\leq & \lim_{\gamma \rightarrow + \infty} \frac{1}{\gamma^{2-h}} M_{2}(\chi).
\end{eqnarray*}
\end{rmk}

Let $\psi$ denotes the class of kernels satisfying the assumptions (K1)-(K2). Then, for $\chi \in \psi$ and $w>0,$ the bivariate Kantorovich exponential sampling series is defined as
\begin{equation} \label{main}
(I^{\chi}_{w} f)(x,y) := \sum_{k= - \infty}^{+\infty} \sum_{j= - \infty}^{+\infty} \chi(e^{-k}x^w, e^{-j}y^w) w^{2}  \int_{\frac{k}{w}}^{\frac{k+1}{w}} \int_{\frac{j}{w}}^{\frac{j+1}{w}} f(e^u,e^v)\  du \ dv
\end{equation}
where $f : \mathbb{R}^{2}_{+} \rightarrow \mathbb{R} $ is a locally integrable function such that the above series is convergent $\forall (x,y) \in \mathbb{R}^{2}_{+}.$ In view of condition (K2), the series (\ref{main}) is well-defined for the class of bounded functions on $\mathbb{R}^{2}_{+}.$

\section{Convergence Results}
In this section, we analyse the pointwise and uniform convergence results and a Voronovskaya type asymptotic formula for the sampling series $(I_{w}^{\chi}f)_{w>0}$ exploiting the bivariate Mellin Taylor's formula.

\begin{thm}\label{theorem1}
Let $\chi \in \psi$ and $ f\in  L^{\infty}(\mathbb{R}^{2}_{+}).$ Then the series $(I^{\chi}_{w} f)(x,y)$ converges to $f(x,y)$ at every point $(x,y) \in \mathbb{R}^{2}_{+},$ the point of continuity of $f.$ Further, for $f \in \mathcal{C}(\mathbb{R}^{2}_{+})$ we have
$$ \lim_{w \rightarrow \infty} \|I_{w}^{\chi}f - f \|_{\infty} = 0.$$
\end{thm}

\noindent\bf{Proof.}\rm \
In view of condition (K1), we write
\begin{eqnarray*}
|I_{w}^{\chi}f(x,y)-f(x,y)|
&\leq& \sum_{k= - \infty}^{+\infty} \sum_{j= - \infty}^{+\infty}\bigg|\chi(e^{-k} x^{w},e^{-j} y^{w})\bigg|w^{2} \int_{\frac{k}{w}}^{\frac{k+1}{w}}\int_{\frac{j}{w}}^{\frac{j+1}{w}}|f(e^{u},e^{v}) - f(x,y)|\ du\  dv \\
&\leq& \left( \sum \sum_{(\frac{k}{w},\frac{j}{w}) \in B_{\frac{\delta}{2}}(\log x,\log y)} + \sum \sum_{(\frac{k}{w},\frac{j}{w}) \notin B_{\frac{\delta}{2}}(\log x,\log y)} \right) \big| \chi(e^{-k} x^{w},e^{-j} y^{w})\big|\\
\\&& \ w^{2} \int_{\frac{k}{w}}^{\frac{k+1}{w}}\int_{\frac{j}{w}}^{\frac{j+1}{w}} |f(e^{u},e^{v}) - f(x,y)| \ dudv
:=I_{1}+I_{2}.
\end{eqnarray*}
As $f \in \mathcal{C}(\mathbb{R}^{2}_{+}),$ for any $ \epsilon >0$ there exists $\delta>0$ such that $ |f(e^{u},e^{u}) - f(x,y)|< \epsilon,$ whenever $\sqrt{(u- \log x)^{2}+(v- \log y)^{2}} < \delta.$
Let $w^{'}$ be fixed such that $ \frac{1}{w}< \frac{\delta}{2} $ for every $ w > w^{'}.$ Now for  $ u \in\big[ \frac{k}{w}, \frac{k+1}{w} \big]$ and $ v \in\big[ \frac{j}{w}, \frac{j+1}{w} \big]$ and $ w > w^{'},$ we have
$$ |u- \log x| \leq \Big | u - \frac{k}{w} \Big|+\Big |\frac{k}{w} - \log x \Big| < \delta,$$
and $$ |v- \log y| \leq \Big | v - \frac{j}{w} \Big|+\Big |\frac{j}{w} - \log y \Big| < \delta,$$
whenever $ \big|\frac{k}{w} - \log x\big|< \frac{\delta}{2}$ and $ \big|\frac{j}{w} - \log y\big|< \frac{\delta}{2}.$
This gives $|I_{1}| < \epsilon M_{0}(\chi).$ Subsequently,  $I_{2}$ can be estimated as
$$ |I_{2}| \leq \ 2 \|f\|_{\infty} \sum \sum_{(\frac{k}{w},\frac{j}{w}) \notin B_{\frac{\delta}{2}}(\log x,\log y)} \big| \chi(e^{-k} x^{w},e^{-j} y^{w})\big|.$$
In view of Remark 2.2, we deduce that $|I_{2}|\leq 2 \epsilon \|f\|_{\infty}.$ Combining the estimates $I_{1}-I_{2}$ we get the desired result.

\begin{thm}\label{theorem2}
Let $ f \in C^{(2)}(\mathbb{R}^{2}_{+})$ locally at $(x,y)$ and $\chi \in \psi$ be the kernel such that $m_{1,0}(\chi)=m_{0,1}(\chi)=0.$ Then we have
$$ \lim_{w\rightarrow \infty} w \big[ (I_{w}^{\chi}f)(x,y) - f(x,y) \big] = \frac{1}{2} \left[ \theta_{x} f(x,y)+\theta_{y} f(x,y) \right].$$
 \end{thm}

\noindent\bf{Proof.}\rm \ Since $f \in C^{(2)}(\mathbb{R}^{+}),$ using the Taylor's formula in terms of Mellin derivatives (see \cite{butzer3,bardaro7}) upto second order term, we can write
\begin{eqnarray*}
f(e^{u},e^{v})&=& f(x,y)+ \theta_{x} f(x,y) (u- \log x)+\theta_{y} f(x,y) (v- \log y) \\&& + \frac{1}{2!}\Big(\theta_{x} (u- \log x)+\theta_{y} (v- \log y)\Big)^{2}f(x,y) + h \Big(\frac{e^{u}}{x},\frac{e^{v}}{y}\Big) \Big( (u- \log x)^{2}+ (v- \log y)^{2}\Big),
\end{eqnarray*}
 where $h$ is a bounded function such that $\displaystyle \lim_{(s,t) \rightarrow (1,1)} h(s,t)=0.$ In view of (\ref{main}), we obtain
\begin{eqnarray*}
(I_{w}^{\chi}f)(x,y) - f(x,y) &=& \sum_{k= - \infty}^{+\infty}\sum_{j= - \infty}^{+\infty} \chi(e^{-k} x^{w},e^{-j} y^{w})\  w^{2} \int_{\frac{k}{w}}^{\frac{k+1}{w}} \int_{\frac{j}{w}}^{\frac{j+1}{w}} \Big [ \theta_{x} f(x,y) (u- \log x)\\
\\&&+ \theta_{y} f(x,y) (v- \log y) + \frac{1}{2!}\Big(\theta_{x} (u- \log x) +\theta_{y} (v- \log y)\Big)^{2}f(x,y) \\&&
+ h \Big(\frac{e^{u}}{x},\frac{e^{v}}{y}\Big) \Big( (u- \log x)^{2}+ (v- \log y)^{2}\Big) \Big] \ du\ dv
:= I_{1}+I_{2}+ I_{3}.
\end{eqnarray*}
It is easy to see that
\begin{eqnarray*}
I_{1}&=& \sum_{k= - \infty}^{+\infty}\sum_{j= - \infty}^{+\infty} \chi(e^{-k} x^{w},e^{-j} y^{w})\  w^{2} \int_{\frac{k}{w}}^{\frac{k+1}{w}} \int_{\frac{j}{w}}^{\frac{j+1}{w}}\Big[\theta_{x} f(x,y) (u- \log x)
 + \theta_{y} f(x,y) (v- \log y) \Big]\ dudv \\
&=& \frac{1}{2w} \left[ \theta_{x} f(x,y)+\theta_{y} f(x,y) \right] .
\end{eqnarray*}
Similarly, we obtain
\begin{eqnarray*}
I_{2} &=& \frac{1}{6w^{2}} \left[2\theta_{x}^{2} f(x,y)\left( 1+3m_{2,0}(\chi)\right)+3\theta_{xy}^{2} f(x,y)\left( 1+4m_{1,1}(\chi)\right)+2\theta_{y}^{2} f(x,y)\left( 1+3m_{0,2}(\chi)\right) \right].
\end{eqnarray*}
In order to estimate $I_{3},$ let $\epsilon > 0$ be fixed then there exists $\delta > 0$ such that $|h(s,t)| < \epsilon $ whenever $|s-1|< \delta$ and $|t-1|< \delta.$ Moreover, let $w^{'}$ be fixed in such a way that $ \frac{1}{w}< \frac{\delta}{2} $ for every $ w > w^{'}.$ We write $I_{3}$ as
\begin{eqnarray*}
|I_{3}|&\leq&\sum \sum_{(\frac{k}{w},\frac{j}{w}) \in B_{\frac{\delta}{2}}(\log x,\log y)} \big |\chi(e^{-k} x^{w},e^{-j} y^{w})\big | \ \ \Bigg| w^{2} \int_{\frac{k}{w}}^{\frac{k+1}{w}} \int_{\frac{j}{w}}^{\frac{j+1}{w}}\left[ h \Big(\frac{e^{u}}{x},\frac{e^{v}}{y}\Big) (u-\log x)^{2}+(v-\log y)^{2} \right] \ du\ dv \Bigg|
\\&& + \sum \sum_{(\frac{k}{w},\frac{j}{w}) \notin B_{\frac{\delta}{2}}(\log x,\log y)} \big|\chi(e^{-k} x^{w},e^{-j} y^{w})\big| \ \ \Bigg| w^{2} \int_{\frac{k}{w}}^{\frac{k+1}{w}} \int_{\frac{j}{w}}^{\frac{j+1}{w}}\left[ h \Big(\frac{e^{u}}{x},\frac{e^{v}}{y}\Big) (u-\log x)^{2}+(v-\log y)^{2} \right] \ du\ dv \Bigg| \\
&:=& I_{3}^{'}+I_{3}^{''}.
\end{eqnarray*}
Using the argument of Theorem \ref{theorem1}, we deduce that
$$ |w I_{3}^{'}| \leq \frac{\epsilon}{6 w}\left( 7 M_{0}(\chi)+24 M_{2}(\chi)+24 M_{1}(\chi) \right).$$
For $I_{3}^{''},$ using Remark 2.2 and the fact that $h$ is bounded, we obtain
\begin{eqnarray*}
|I_{3}^{''}| &\leq & \|h\|_{\infty} \ \sum \sum_{(\frac{k}{w},\frac{j}{w}) \notin B_{\frac{\delta}{2}}(\log x,\log y)} \big |\chi(e^{-k} x^{w},e^{-j} y^{w})|\  w^{2} \int_{\frac{k}{w}}^{\frac{k+1}{w}}\int_{\frac{j}{w}}^{\frac{j+1}{w}} \left((u-\log x)^2 + (v-\log y)^2 \right) \ du\ dv  \\
&\leq & \frac{55 \epsilon \|h\|_{\infty}}{6 w^2} .
\end{eqnarray*}
Hence, we obtain $\displaystyle |w I_{3}^{''}| \leq \frac{55 \epsilon \|h\|_{\infty}}{6 w}.$ Combining the estimates of $I_{1}-I_{3},$ we get required result.

\begin{rmk}
It is important to note that the similar asymptotic formula can be obtained for $f \in C^{(1)}(\mathbb{R}^{2}_{+})$. But, the assumption $f \in C^{(2)}(\mathbb{R}^{2}_{+})$ is considered to show that the family (\ref{main}) converges linearly even if the higher order partial derivatives exist on $\mathbb{R}^{2}_{+}.$
\end{rmk}

\begin{rmk}
The condition that $f$ is bounded on $ \mathbb{R}^{2}_{+}$ in Theorem 4.2 can be relaxed by assuming that there are arbitrary constants $a,b$ such that
$ |f(x,y)| \leq a + b|\log^{2}x + \log^{2}y | , \hspace{0.25cm} \forall x,y \in \mathbb{R}^{2}_{+}.$
\end{rmk} First we show that the series (\ref{main}) is well defined for such $f.$ Indeed
\begin{eqnarray*}
|I_{w}^{\chi}f)(x)| &\leq &  \sum_{k= - \infty}^{+\infty}\sum_{j= - \infty}^{+\infty}  |\chi(e^{-k} x^{w},e^{-j} y^{w})| \  w^{2} \int_{\frac{k}{w}}^{\frac{k+1}{w}} \int_{\frac{j}{w}}^{\frac{j+1}{w}}|f(e^{u},e^{v})| \ dudv \\
& \leq & \sum_{k= - \infty}^{+\infty} \sum_{j= - \infty}^{+\infty}|\chi(e^{-k} x^{w},e^{-j} y^{w})| \  w^{2} \int_{\frac{k}{w}}^{\frac{k+1}{w}} (a + b| u^{2}+v^{2}|) \ dudv \\
& \leq & M_{0}(\chi) \left[ a+b \left(|\log^{2} x| + |\log^{2}y| \right)+ \frac{2b}{3w^{2}} \right]+ \frac{b}{w^{2}} \left[ M_{1,0}(\chi)+M_{0,1}(\chi)+M_{2,0}(\chi)+M_{0,2}(\chi) \right ].
\end{eqnarray*}
This shows that the series $(I_{w}^{\chi}f)_{w>0}$ is absolutely convergent in $\mathbb{R}^{2}_{+}.$ For any fixed $(x,y) \in \mathbb{R}^{2}_{+},$ we define
\begin{eqnarray*}
P_{2}(u,v)&:=& f(x,y)+ (\theta_{x} f)(x,y) (u-\log x)+(\theta_{y} f)(x,y) (v-\log y)\\
& +& \frac{1}{2!}\left((\theta^{2}_{x} f)(x,y)(u-\log x)^{2}+2(\theta^{2}_{xy} f)(x,y)(u-\log x)(v-\log y)+(\theta^{2}_{y} f)(x,y)(v-\log y)^{2}\right).
\end{eqnarray*}
From the bivariate Taylor's formula in terms of Mellin derivatives upto second order term, we can write as
$$ h \bigg(\frac{e^{u}}{x},\frac{e^{u}}{y}\bigg) = \frac{f(e^{u},e^{v}) - P_{2}(u,v)}{(u-\log x)^{2}+(v-\log y)^{2}},$$
where $h$ is a function such that $\displaystyle \lim_{\stackrel{(u- \log x) \rightarrow 0}{(v- \log y) \rightarrow 0}} h \bigg(\frac{e^{u}}{x},\frac{e^{v}}{y}\bigg)=0.$ This implies that $h$ is bounded in $\delta-$ neighbourhood of $(\log x, \log y),$ i.e $(u,v) \in B_{\delta}(\log x, \log y).$ Now for $(u,v) \notin B_{\delta}(\log x, \log y),$ we have
\begin{eqnarray*}
|h(e^{u} x^{-1},e^{v} y^{-1})| & \leq & \frac{|f(e^{u},e^{v})|}{|u - \log x|^{2}+|v - \log y|^{2}} + \frac{|P_{2}(u,v)|}{|u - \log x|^{2}+|v - \log y|^{2}} \\
& \leq & \frac{a+ b |u^{2}+v^{2}|}{|u - \log x|^{2}+|v - \log y|^{2}} + \frac{|P_{2}(u,v)}{|u - \log x|^{2}+|v - \log y|^{2}} \ .
\end{eqnarray*}
This show that  $h(.,.)$ is also bounded for $(u,v) \notin B_{\delta}(\log x, \log y).$ This concludes that $ h$ is bounded on $\mathbb{R}^{2}_{+}.$ Now we can proceed in the similar manner as in the proof of Theorem 4.2 to get the same asymptotic formula.\\

For $f \in \mathcal{C}(\mathbb{R}^{2}_{+}),$ the logarithmic modulus of continuity is defined as
$$ \omega(f,\delta):= \sup \{|f(x,y)-f(u,v)|: \  |\log x-\log u| \leq \delta_{1},|\log y-\log v| \leq \delta_{2}\  \ \delta_{1},\delta_{2} \in \mathbb{R}^{+}\} .$$
For every $\delta_{1} > 0,$ $\delta_{2} > 0$ and $(x,y),(u,v) \in \mathbb{R}^{2}_{+},$ \ $\omega$ \ satisfies the following properties:
\begin{itemize}
\item[a)] $\omega(f, \delta_{1},\delta_{2}) \rightarrow 0$ as $\delta_{1} \rightarrow 0$ and $\delta_{2} \rightarrow 0 .$
\item[b)] $|f(x,y) - f(u,v)| \leq \omega(f,\delta_{1},\delta_{2}) \left( 1+ \frac{|\log x - \log u|}{\delta_{1}} \right)\left( 1+ \frac{|\log y - \log v|}{\delta_{2}} \right).$
\end{itemize}
Further properties of modulus of continuity can be found in \cite{mamedeo,bardaro9}. Now we obtain a quantitative estimate of the convergence of series (\ref{main}) for $ f \in \mathcal{C}(\mathbb{R}^{2}_{+}).$

\begin{thm}\label{t3}
Let $ f \in \mathcal{C}(\mathbb{R}^{2}_{+}).$ Then for any $(x,y)\in\mathbb{R}^{2}_{+}),$ we have
\begin{eqnarray*}
|(I_{w}^{\chi}f)(x,y) - f(x,y)| &\leq &  \omega(f,\delta_{1},\delta_{2}) \left[  M_{0}(\chi) \left(  1+ \frac{1}{2 \delta_{1}w}+ \frac{1}{2 \delta_{2}w}+ \frac{1}{4 \delta_{1} \delta_{1} w^{2}}\right) + \frac{M_{1,0}(\chi)}{\delta_{1}w} \left(1+ \frac{1}{2w \delta_{1}}  \right) \right] \\&&
+ \ \omega(f,\delta_{1},\delta_{2}) \left[ \frac{M_{0,1}(\chi)}{\delta_{2}w} \left(1+ \frac{1}{2w \delta_{2}}  \right)+  \frac{M_{1,1}(\chi)}{\delta_{1} \delta_{2} w^{2}} \right],
\end{eqnarray*}
for any $w>0$ and $\delta_{1}>0$,$\delta_{2}>0.$
\end{thm}

\noindent\bf{Proof.}\rm \ Using property (b) of $\omega$, we can write
\begin{eqnarray*}
|(I_{w}^{\chi}f)(x,y) - f(x,y)| & \leq & \omega(f,\delta_{1},\delta_{2})\sum_{k= - \infty}^{+\infty} \sum_{j= - \infty}^{+\infty} |\chi(e^{-k} x^{w},e^{-j} y^{w})| \ w^{2}
\\&& \int_{\frac{k}{w}}^{\frac{k+1}{w}}\int_{\frac{j}{w}}^{\frac{j+1}{w}} \left( 1+ \frac{\left|u -  \log x \right|}{\delta_{1}} \right)\left( 1+ \frac{\left|v -  \log y \right|}{\delta_{2}} \right)\ du \ dv \\
& \leq & \omega(f,\delta_{1},\delta_{2})\sum_{k= - \infty}^{+\infty} \sum_{j= - \infty}^{+\infty} |\chi(e^{-k} x^{w},e^{-j} y^{w})| w^{2}
\\&& \int_{\frac{k}{w}}^{\frac{k+1}{w}}\int_{\frac{j}{w}}^{\frac{j+1}{w}}\left( 1+\frac{\left|u -  \log x \right|}{\delta_{1}}+\frac{\left|v -  \log y \right|}{\delta_{2}} +\frac{\left|u -  \log x \right|}{\delta_{1}}\frac{\left|v -  \log y \right|}{\delta_{2}} \right)du\ dv\\
&:=& J_{1}+J_{2}+J_{3}+J_{4}.
\end{eqnarray*}
It is easy to see that $J_{1}= M_{0}(\chi) \omega(f,\delta_{1},\delta_{2}).$ Now we estimate $J_{2}.$
\begin{eqnarray*}
J_{2} &\leq & \frac{\omega(f,\delta_{1},\delta_{2})}{\delta_{1}} \sum_{k= - \infty}^{+\infty} \sum_{j= - \infty}^{+\infty} |\chi(e^{-k} x^{w},e^{-j} y^{w})|\  w^{2} \int_{\frac{k}{w}}^{\frac{k+1}{w}}\int_{\frac{j}{w}}^{\frac{j+1}{w}}  |u-\log x| \ du\ dv \\
&\leq & \frac{\omega(f,\delta_{1},\delta_{2})}{2 \delta_{1}w} \left[ M_{0}(\chi)+ 2 M_{1,0}(\chi) \right].
\end{eqnarray*}
Similarly, we obtain $ \displaystyle J_{3} \leq \frac{\omega(f,\delta_{1},\delta_{2})}{2 \delta_{2}w} \left[ M_{0}(\chi)+ 2 M_{0,1}(\chi) \right ].$
Finally, we evaluate $J_{4}.$
\begin{eqnarray*}
J_{4} &\leq & \frac{\omega(f,\delta_{1},\delta_{2})}{\delta_{1} \delta_{2}} \sum_{k= - \infty}^{+\infty} \sum_{j= - \infty}^{+\infty} |\chi(e^{-k} x^{w},e^{-j} y^{w})| \ w^{2} \int_{\frac{k}{w}}^{\frac{k+1}{w}}\int_{\frac{j}{w}}^{\frac{j+1}{w}}  |u-\log x| |v-\log y|\ du\ dv \\
&\leq & \frac{\omega(f,\delta_{1},\delta_{2})}{4 \delta_{1} \delta_{2} w^{2}} \left[ M_{0}(\chi)+2 M_{1,0}(\chi)+2 M_{0,1}(\chi)+4 M_{1,1}(\chi) \right].
\end{eqnarray*}
On combining the estimates $J_{1}-J_{4},$ we get the desired estimate.

\begin{rmk}
For fixed $w>0,$ if we put $\delta_{1}=\delta_{2}= \frac{1}{w}$ then the estimate in Theorem 3.3 is as follows:
$$ |(I_{w}^{\chi}f)(x,y) - f(x,y)| \leq \frac{\kappa}{4} \ \omega \left(f, \frac{1}{w} \right)\ , \hspace{0.5cm} \kappa := \left( 9 M_{0}(\chi)+ 6 M_{1,0}(\chi)+6 M_{0,1}(\chi)+ 4 M_{1,1}(\chi) \right).$$
\end{rmk}

\section{GBS-Bivariate Kantorovich Sampling Series}
 B\"{o}gel pioneered the notion of $B$-continuous and $B$-differentiable functions and furnished significant results in \cite{BOG1,BOG2}. Badea and Cottin \cite{badea,BAD2} established the well-known Korovkin theorem for B-continuous functions. Dobrescu and Matei \cite{DOB} showed that $B$-continuous function can be approximated uniformly by GBS operators associated to Bernstein polynomials on bounded domain. The study of approximation behaviour of GBS operators of Bernstein-Stancu polynomials was developed in \cite{MIC}. Agrawal et al. \cite{PNA6} derived the order of convergence for Lupa\c{s}-Durrmeyer type GBS operators on the basis of P\'{o}lya distribution. Since then, various authors have contributed significantly in this direction, see eg \cite{PNA3,BAR,KAJ,POP,own1,ruchi,acar}.

For $X, Y \subseteq \mathbb{R}_{+},$ we call a function $f:X\times Y \rightarrow \mathbb{R}\ $ as Mellin $B$-continuous (Mellin B\"{o}gel Continuous) at $(u,v)\in X \times Y $ if
$$\lim_{(x,y)\rightarrow (e^u,e^v)}\Delta_{(x,y)}f[e^u,e^v;x,y]=0,$$
where $\Delta_{(x,y)}f[e^u,e^v;x,y]=f(x,y)-f(x,e^v)-f(e^u,y)+f(e^u,e^v).$ The function $f:X\times Y \rightarrow \mathbb{R}$ is Mellin $B$-bounded if there exists $\lambda>0$ such that $\mid\Delta_{(x,y)}f[e^u,e^v;x,y]\mid \ \leq \lambda$ for every $(x,y), (u,v)\in X\times Y .$
Throughout this paper, we denote $\mathcal{B}_{b}(\mathbb{R}^{2}_{+})$ and $\mathcal{C}_{b}(\mathbb{R}^{2}_{+})$ be the space of all Mellin $B$-bounded and Mellin $B$-continuous functions on $\mathbb{R}^{2}_{+}$ with the usual sup-norm $ \|. \|_{\infty}$ respectively.

Now for any $f\in \mathcal{C}_{b}(\mathbb{R}{^2}_{+}),$ the GBS of bivariate Kantorovich exponential sampling series is defined as
$$ \tilde{I}_{w}^{\chi}(f;x,y) :=  I_{w}^{\chi}(f(x,v)+f(u,y)-f(u,v); x,y)$$
for all $ (x,y), \ (u,v) \in \mathbb{R}^{2}_{+}.$

We estimate the order of convergence for the family of operators $\tilde{I}_w^{\chi}(f;x,y)$ using B\"{o}gel frame of modulus of continuity for $f\in \mathcal{C}_{b}(\mathbb{R}^{2}_{+}).$ For any $\delta_1,\delta_2 > 0 ,$ the mixed modulus of smoothness in the Mellin sense is defined as
$$\omega_{B}(f;\delta_1,\delta_2):=\sup\{\mid\Delta_{(x,y)}f[s,t;x,y]\mid \ : \ \mid \log s - \log x \mid <\delta_1, \mid \log t - \log y \mid <\delta_2\}$$ for every $(x,y), (u,v)\in\mathbb{R}^{2}_{+}.$
\begin{thm}\label{t3}
Let $f\in \mathcal{C}_{b}(\mathbb{R}^{2}_{+}).$ Then the following estimate holds
$$ \mid \tilde{I}^{\chi}_{w}(f;x,y)-f(x,y)\mid \leq \left(1+\frac{A_{1}}{\delta_{1}}+\frac{A_{2}}{\delta_{2}}+\frac{A_{3}}{\delta_{1}\delta_{2}}\right)\omega_{B}\left(f;\delta_{1},\delta_{2}\right)$$
where $A_{1}=\frac{1}{2w}\left(M_{0,0}+2M_{1,0}\right), A_{2}=\frac{1}{2w}\left(M_{0,0}+2M_{0,1}\right), A_{3}=\frac{1}{4w^2}\left(M_{0,0}+2M_{1,0}+2M_{0,1}+4M_{1,1}\right).$
\end{thm}

\noindent\bf{Proof.}\rm \ Using the property $\omega_{B}\left(f;\alpha \delta_{1},\beta \delta_{2}\right)\leq(1+\alpha)(1+\beta)\omega_{B}\left(f;\delta_{1},\delta_{2}\right)$ for $\alpha,\beta >0,$ we write
\begin{eqnarray*}
\mid \Delta_{(x,y)}f[e^u,e^v;x,y]\mid &\leq&\omega_{B}\left(f;\mid u-\log x\mid,\mid v-\log y\mid\right) \\
&\leq & \left(1+\frac{\mid u-\log x\mid}{\delta_{1}}\right)\left(1+\frac{\mid v-\log y\mid}{\delta_{2}}\right)\omega_{B}\left(f;\delta_{1},\delta_{2}\right)
\end{eqnarray*}
for any $(x,y),(u,v) \in \mathbb{R}^{2}_{+}$ and $\delta_{1},\delta_{2} >0.$ Now, applying the series $(\tilde{I}^{\chi}_{w})$ on $\Delta_{(x,y)}f[e^u,e^v;x,y],$ we obtain
\begin{eqnarray}\label{condition}
\tilde{I}^{\chi}_{w}(f;x,y)= f(x,y)-I_{w}^{\chi}(\Delta_{(x,y)}f[e^u,e^v;x,y]).
\end{eqnarray}
From (\ref{condition}), we have
\begin{eqnarray*}
\mid(\tilde{I}^{\chi}_{w})(f;x,y)-f(x,y) \mid &\leq & \Bigg[ 1+\frac{I_{w}^{\chi}(\mid u-\log x\mid;x,y)}{\delta_{1}}+\frac{I_{w}^{\chi}(\mid v-\log y\mid;x,y)}{\delta_{2}}\\
\\&& \ \ \ \  +\ \frac{I_{w}^{\chi}(\mid u-\log x\mid \mid v-\log y\mid;x,y)}{\delta_{1}\delta_{2}} \Bigg] \omega_{B}\left(f;\delta_{1},\delta_{2}\right).
\end{eqnarray*}
Using the definition (\ref{main}), we obtain
\begin{eqnarray*}
I_{w}^{\chi}(\mid u-\log x\mid;x,y)&=& \sum_{k= - \infty}^{+\infty} \sum_{j= - \infty}^{+\infty} \chi(e^{-k} x^{w},e^{-j} y^{w}) w^{2}
 \int_{\frac{k}{w}}^{\frac{k+1}{w}}\int_{\frac{j}{w}}^{\frac{j+1}{w}}\mid u-\log x\mid du\ dv \\
 &=&\frac{1}{2w}\sum_{k= - \infty}^{+\infty} \sum_{j= - \infty}^{+\infty} \chi(e^{-k} x^{w},e^{-j} y^{w})[1+2\left(w\log x-k\right)] \\
&=&\frac{1}{2w}\left( M_{0,0}+2M_{1,0}\right).
\end{eqnarray*}
Similarly, we get $I_{w}^{\chi}(\mid v-\log y\mid;x,y)=\frac{1}{2w}\left( M_{0,0}+2M_{0,1}\right).$ Finally, we have
\begin{eqnarray*}
I_{w}^{\chi}(\mid u-\log x\mid\mid v-\log y\mid;x,y)&=& \sum_{k= - \infty}^{+\infty} \sum_{j=-\infty}^{+\infty} \chi(e^{-k} x^{w},e^{-j} y^{w}) w^{2}\int_{\frac{k}{w}}^{\frac{k+1}{w}}\int_{\frac{j}{w}}^{\frac{j+1}{w}}\mid u-\log x\mid \mid v-\log y\mid du\ dv \\
 &=&\frac{1}{2w}\sum_{k= - \infty}^{+\infty} \sum_{j= - \infty}^{+\infty} \chi(e^{-k} x^{w},e^{-j} y^{w})[1+2\left(w\log y-k\right)][1+2\left(w\log x-k\right)] \\
 &=&\frac{1}{4w^{2}}\left( M_{0,0}+2M_{1,0}+2M_{0,1}+4M_{1,1}\right).
\end{eqnarray*}
On substituting these estimates, we get the desired result.\\

Now we define a Mellin B-differential (Mellin B\"{o}gel differential) function. A function $f:X\times Y(\subseteq \mathbb{R}^{2}_{+}) \rightarrow \mathbb{R}$ is said to be Mellin B-differential at $(u,v),$ if the following limit exist finitely
$$\lim_{(x,y)\longrightarrow (e^u,e^v)} \frac{\Delta_{(x,y)}f[e^u,e^v;x,y]}{(u- \log x)(v- \log y)}.$$
The Mellin B-differential of $f$ at any point $(u,v)$ is represented as $\theta_{B}(f;u,v).$ Moreover, we denote $\mathcal{D}_{b}(\mathbb{R}^{2}_{+})$ as the space of all Mellin B-differentiable function on $ \mathbb{R}^{2}_{+}.$

\begin{thm} Let $f\in \mathcal{D}_{b}(\mathbb{R}^{2}_{+})$ and $\theta_{B}f \in \mathcal{B}(\mathbb{R}^{2}_{+}).$ Then for each $(x,y)\in\mathbb{R}^{2}_{+},$ we have
\begin{eqnarray*}
 \mid(\tilde{I}^{\chi}_{w})(f;x,y)-f(x,y) \mid &\leq& E_{1} \left( 3\|\theta_{B} f \|_{\infty}+\omega_{B}\left(f;\delta_{1},\delta_{2}\right)\right)+\left(\frac{E_{2}}{\delta_{1}}+\frac{E_{3}}{\delta_{2}}+\frac{E_{4}}{\delta_{1}\delta_{2}}\right)\omega_{B}\left(\theta_{B}f;\delta_{1},\delta_{2}\right)
 \end{eqnarray*}
where $E_{1}=\frac{1}{4w^{2}}\left( M_{0,0}+2M_{1,0}+2M_{0,1}+4M_{1,1}\right),E_{2}=\frac{1}{6w^{3}}\left( M_{0,0}+3M_{2,0}+3M_{1,0}+2M_{0,1}+6M_{2,1}+6M_{1,1}\right),\\E_{3}=\frac{1}{6w^{3}}\left(M_{0,0}+3M_{0,2}+3M_{0,1}+2M_{1,0}+6M_{1,2}+6M_{1,1}\right),
\\E_{4}=\frac{1}{9w^{4}}\left(M_{0,0}+3M_{2,0}+3M_{0,2}+3M_{1,0}+3M_{0,1}+
9M_{2,2}+9M_{1,2}+9M_{2,1}+9M_{1,1}\right).$
 \end{thm}

\noindent\bf{Proof.}\rm \ Since $f\in \mathcal{D}_{b}(\mathbb{R}^{2}_{+}),$ we have
$\Delta_{(x,y)}f[e^u,e^v;x,y]=(u-\log x)(v-\log y) \theta_{B}f(p,q),$
where $p\in(\log x,u)$ and $q \in (\log y,v).$ Using the definition of $\Delta_{(x,y)}f[e^u,e^v;x,y],$ we write
\begin{eqnarray*}
\theta_{B}f(p,q)&=&\Delta_{(x,y)}\theta_{B}f(p,q)+\theta_{B}f(p,y)+\theta_{B}f(x,q)-\theta_{B}f(x,y).
\end{eqnarray*}
Using the above equality and the fact that $\theta_{B}f \in \mathcal{B}(\mathbb{R}^{2}_{+}),$ we have
\begin{eqnarray*}
\mid I_{w}^{\chi} \left(\Delta_{(x,y)}f[e^u,e^v;x,y];x,y \right)\mid &=& \mid I_{w}^{\chi} \left((u-\log x)(v-\log y)\theta_{B}f(p,q);x,y \right)\mid\\
&\leq& I_{w}^{\chi} \left(\mid u-\log x\mid\mid v-\log y\mid \mid \Delta_{x,y}f(p,q)\mid;x,y \right)
\\&&+I_{w}^{\chi} \left(\mid u-\log x\mid\mid v-\log y\mid \left(\mid \theta_{B}f(p,y)\mid+\mid \theta_{B}f(x,q)\mid+\mid \theta_{B}f(x,y)\mid\right);x,y \right)\\
&\leq& I_{w}^{\chi} \left(\mid u-\log x\mid\mid v-\log y\mid \omega_{B}\left(\theta_{B}f;\mid p-\log x \mid,\mid q-\log y \mid\right)\right)
\\&&+3\|\theta_{B}f\|_{\infty}I_{w}^{\chi}\left(\mid u-\log x\mid\mid v-\log y\mid \right).
\end{eqnarray*}
Using the monotonicity of mixed modulus of smoothness $\omega_{B},$ we write
\begin{eqnarray*}
\mid(\tilde{I}^{\chi}_{w})(f;x,y)-f(x,y) \mid &\leq& 3\|\theta_{B}f\|_{\infty} \ I_{w}^{\chi}\left(\mid u-\log x\mid\mid v-\log y \mid;x,y\right)
 \\&&+\Bigg[\frac{I_{w}^{\chi}\left(\mid u-\log x\mid^{2}\mid v-\log y\mid;x,y\right)}{\delta_{1}}+\frac{I_{w}^{\chi}\left(\mid u-\log x\mid\mid v-\log y\mid^{2};x,y\right)}{\delta_{2}}
 \\&&+\frac{I_{w}^{\chi}\left(\mid u-\log x\mid^{2}\mid v-\log y\mid^{2};x,y\right)}{\delta_{1}\delta_{2}}\Bigg ] \omega_{B}\left(\theta_{B}f;\delta_{1},\delta_{2}\right).
 \end{eqnarray*}
From the definition (\ref{main}), we obtain
\begin{eqnarray*}
 I_{w}^{\chi}\left(\mid u-\log x\mid\mid v-\log y\mid;x,y\right)&=&\frac{1}{4w^{2}}\Big[ M_{0,0}+2M_{1,0}+2M_{0,1}+4M_{1,1}\Big]. \\
 I_{w}^{\chi}\left(\mid u-\log x\mid^{2}\mid v-\log y\mid;x,y\right)&=&\frac{1}{6w^{3}}\Big[ M_{0,0}+3M_{2,0}+3M_{1,0}+2M_{0,1}+6M_{2,1}+6M_{1,1}\Big].\\
 I_{w}^{\chi}\left(\mid u-\log x\mid\mid v-\log y\mid^{2};x,y\right)&=&\frac{1}{6w^{3}}\Big[M_{0,0}+3M_{0,2}+3M_{0,1}+2M_{1,0}+6M_{1,2}+6M_{1,1}\Big].\\
I_{w}^{\chi} \left(\mid u-\log x\mid^{2}\mid v-\log y\mid^{2};x,y\right)&=&\frac{1}{9w^{4}}\Big[M_{0,0}+3M_{2,0}+3M_{0,2}+3M_{1,0}+3M_{0,1}+9M_{2,2}
\\&& \hspace{0.9cm} +9M_{1,2}+9M_{2,1}+9M_{1,1}\Big].
 \end{eqnarray*}
Using these estimates, we obtain the required result.\\

Now, we study the degree of approximation for the series by $B$-continuous
functions belonging to the Lipschitz class. The Lipschitz class is given by
\begin{eqnarray*}
Lip_{K}&=&\{f\in \mathcal{C}_{b}(\mathbb{R}^{2}_{+});\mid\Delta_{(x,y)}f[s,t;x,y]\mid\leq K \mid \log s- \log x \mid\  \mid \log t - \log y \mid, \ L \in \mathbb{R}^{+} \}
\end{eqnarray*}
where $(u,v),(x,y) \in \mathbb{R}^{2}_{+}.$
\begin{thm} Let $f \in Lip_{K}$. Then the following holds
\begin{eqnarray*}
\mid \tilde{I}^{\chi}_{w}(f;x,y)-f(x,y)\mid &\leq & \frac{K}{4w^{2}}\left( M_{0,0}+2M_{1,0}+2M_{0,1}+4M_{1,1}\right)
\end{eqnarray*}
where $K$ is the Lipschitz constant.
 \end{thm}

\noindent\bf{Proof.}\rm \ From the definition (\ref{main}), we can write
\begin{eqnarray*}
\tilde{I}^{\chi}_{w}(f;x,y) &=& I_{w}^{\chi}\left( f(x,u)+f(v,y)-f(u,v)\right)\\
&=&I_{w}^{\chi}\left(f(x,y)-\Delta_{(x,y)}f[e^u,e^v;x,y];x,y\right)\\
&=& f(x,y) I_{w}^{\chi}(\textbf{1};x,y)-I_{w}^{\chi}\left(\Delta_{(x,y)}f[e^u,e^v;x,y];x,y\right)
\end{eqnarray*}
Using the estimates of $I_{w}^{\chi} \left( \mid u-\log x\mid \mid v-\log y\mid;x,y\right),$ we obtain
\begin{eqnarray*}
\mid \tilde{I}^{\chi}_{w}(f;x,y)-f(x,y)\mid &\leq&\frac{K}{4w^{2}}\left( M_{0,0}+2M_{1,0}+2M_{0,1}+4M_{1,1}\right).
\end{eqnarray*}
This completes the proof.

\section{Examples of the kernels}
In this section we provide a few examples of the kernel function in the setting of Mellin theory satisfying the assumptions (K1)-(K2). We start with the well-known Mellin-B spline kernels \cite{bardaro9,bardaro7}. For $ x \in \mathbb{R}^{+},$ the $n^{th}$ order Mellin B-spline function are defined as
$$\bar{B}_{n}(x):= \frac{1}{(n-1)!} \sum_{j=0}^{n} (-1)^{j} {n \choose j} \bigg( \frac{n}{2}+\log x-j \bigg)_{+}^{n+1}.$$
$\bar{B}_{n}(x) $ is compactly supported for every $ n \in \mathbb{N}.$ The Mellin transformation of $\bar{B}_{n}$ (see \cite{bardaro7}) is given by
\begin{eqnarray} \label{splinefourier}
\hat{M}[\bar{B}_{n}](c+is) = \bigg( \frac{\sin(\frac{s}{2})}{(\frac{s}{2})} \Bigg)^{n},  \ \ \hspace{0.5cm} s \neq 0.
\end{eqnarray}
Indeed, $\hat{B}_{n}(x)$ satisfies the assumptions (K1)-(K2) for every $n \in \mathbb{N}$ and $x \in \mathbb{R}_{2}^+,$ see \cite{bardaro7,own} .
We can produce the multivariate kernel using univariate kernel functions as follows :
$$ \hat{\chi}^{n}(x_{1},x_{2},...,x_{n}) := \prod_{i=1}^{n} \chi(x_{i})$$
where $ x_{i} \in \mathbb{R}_{+}^{2},$ and $\chi(.)$ is any univariate kernel. For $n=2,$ we construct the bivariate kernel in variables $x$ and $y$ as
$$ \hat{\chi}^{2}(x,y) := \chi(x) \chi(y)\ , \hspace{0.4cm}(x, y) \in \mathbb{R}_{+}^{2}.$$
Exploiting the above construction, we define the bivariate B-spline kernel by
\begin{equation*}
\hat{B}_{2}^{2}(x,y)=\hat{B}_{2}(x)\hat{B}_{2}(y) =
     \begin{cases}

     {(1+\log x)(1+\log y),} &\quad\text{if} \ \ \ \ { e^{-1} < x,y < 1}\\
       {(1-\log x)(1+\log y),} &\quad\text{if} \ \ \ \ { 1 < x <e,\ e^{-1} < y < 1} \\
       {(1+\log x)(1-\log y),} &\quad\text{if} \ \ \ \ { e^{-1} < x < 1,\ 1 < y <e} \\
       {(1-\log x)(1-\log y),} &\quad\text{if} \ \ \ \ { 1 < x,y <e} \\
       {0,} &\quad\text{if} \ \ \ \ { \text{otherwise}.}
   \end{cases}
\end{equation*}
We show the approximation of $f(x,y)= \sin(x^{2}-y^{2}),\ (x,y)\in [0,2]\times[0,2]$ by $(I_{w}^{\chi}f)_{w>0}$ using $\hat{B}_{2}^{2}(x,y)$ as the kernel function (see \textit{Figure-1} and \textit{Table-1}).

\begin{figure}[h]
\centering
{\includegraphics[width=0.8\textwidth]{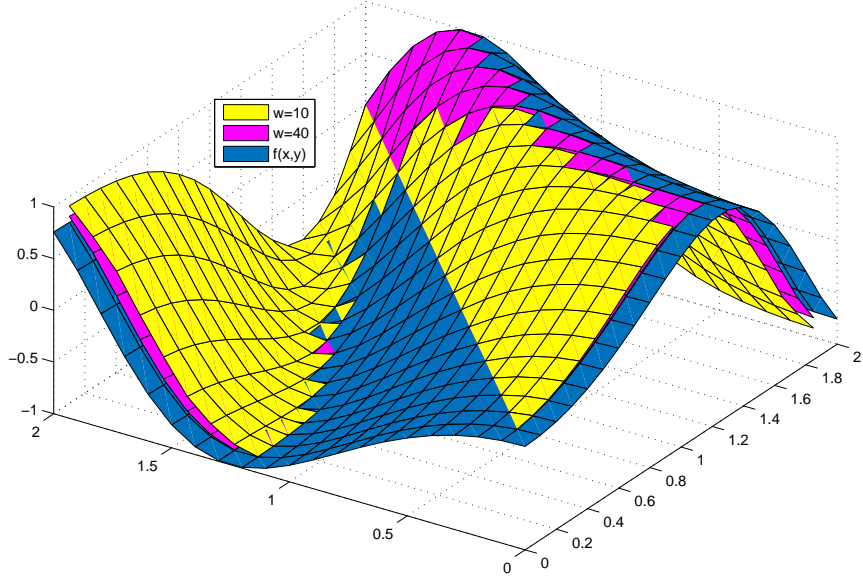}}
\caption{This figure shows the approximation of $f(x,y)$ (Blue) by the series $(I_{w}^{\chi}f)_{w>0}$ for $ w=10$ and $40$ (Yellow and Pink respectively).}
\end{figure}

\begin{tb}\label{table1}\centering
{\it Error estimation (upto $4$ decimal points) in the approximation of $f(x,y)$ by $I_{w}^{\chi}f(x,y)$ for $w=10,40.$}

$  $

\begin{tabular}{|l|l|l|l|}\hline
 $x$&$y$&$ |f(x,y) - I_{10}^{\chi}f(x,y)|$&$|f(x,y)-I_{40}^{\chi}f(x,y)|$\\
 \hline
 $0.2$&$0.1$ & $0.0033$ & $0.0008$\\
  \hline
 $0.6$&$0.5$ & $0.0119$ & $0.0028$ \\
  \hline
 $1.1$ &$0.9$ &$0.0366$ & $0.0092$\\
  \hline
 $1.9$ & $1.8$ & $0.0185$ & $0.0052$  \\
 \hline
             \end{tabular}
   \end{tb}

Next we show the convergence of the function $g(x,y)=y^{2}+\cos(\pi x),\ (x,y)\in [1,4]\times[1,4]$ by the series $(I_{w}^{\chi}g)_{w>0}$ (see \textit{Figure-2} and \textit{Table-2}).
\begin{figure}[h]
\centering
{\includegraphics[width=0.7 \textwidth]{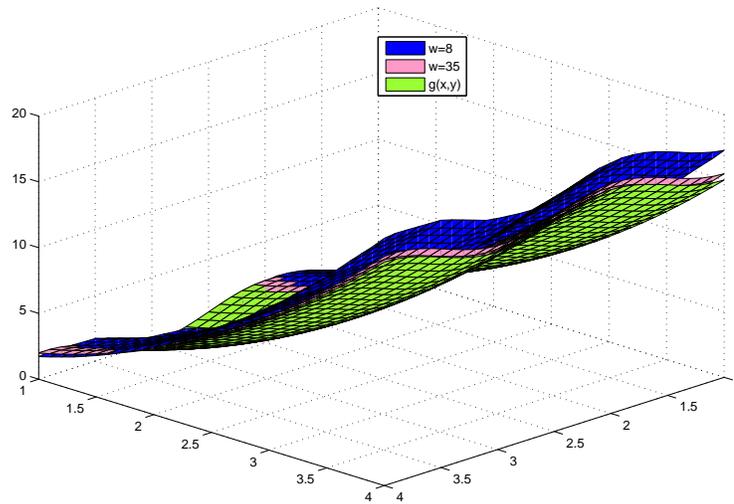}}
\caption{This figure shows the convergence of $g(x,y)$ (Green) by the series $(I_{w}^{\chi} g)_{w>0}$ for $ w= 8$ and $35 $ (Blue and Pink respectively).}
\end{figure}

\begin{tb}\label{table2}\centering
 {\it Error estimation (upto $4$ decimal points) in the convergence of $g(x,y)$ by $I_{w}^{\chi}g(x,y)$ for $w=5,35.$}

$  $

\begin{tabular}{|l|l|l|l|}\hline
 $x$&$y$&$ |g(x,y) - I_{8}^{\chi}g(x,y)|$&$|g(x,y)-I_{35}^{\chi}g(x,y)|$\\
 \hline
 $1.3$&$1.6$ & $0.6079$ & $0.1254$\\
  \hline
 $1.9$&$1.7$ & $0.3958$ & $0.1036$ \\
  \hline
 $2.8$ &$2.4$ &$0.7695$ & $0.1057$\\
  \hline
 $3.6$ &$3.9$ & $2.4188$ & $0.5949$  \\
 \hline
             \end{tabular}
   \end{tb}

\end{document}